\documentclass[runningheads]{llncs}
\usepackage[ruled,vlined,algochapter,linesnumbered, boxed]{algorithm2e}
\usepackage{amsmath,amssymb,mathabx}
\usepackage[T1]{fontenc}
\usepackage{booktabs}
\usepackage{cleveref}
\usepackage{pdflscape}
\usepackage{afterpage}
\usepackage{capt-of}
\usepackage{xspace}
\usepackage{mathtools}
\usepackage{etoolbox}
\usepackage{multirow}
\usepackage{url}
\usepackage[table]{xcolor}

\usepackage{enumitem}
\setlist[enumerate]{itemsep=0mm}

\newcommand{\fracx}[1]{\bar{x}^{(#1)}}

\newcommand{\I}{\mathcal{I}}
\newcommand{\J}{\mathcal{J}}

\usepackage[textsize=scriptsize,textwidth=2.3cm]{todonotes}

\newcommand{\noprint}[1]{}
\usepackage{graphicx}
\begin{document}
\title{A Multi-Reference Relaxation Enforced Neighborhood Search Heuristic in SCIP}
\titlerunning{MRENS Primal Heuristic in SCIP}
\author{
Suresh Bolusani$^*$\inst{1}\orcidID{0000-0002-5735-3443} \and
Gioni Mexi\inst{1}\orcidID{0000-0003-0964-9802}\and
Mathieu Besançon\inst{2}\orcidID{0000-0002-6284-3033} \and
Mark Turner\inst{1}\orcidID{0000-0001-7270-1496}
}
\authorrunning{S. Bolusani et al.}
\institute{Zuse Institute Berlin, Germany \and Grenoble Alpes University, Inria, LIG\\
\email{\{bolusani, mexi, turner\}@zib.de, mathieu.besancon@inria.fr}
}

\let\oldmaketitle\maketitle
\renewcommand{\maketitle}{\oldmaketitle\setcounter{footnote}{0}}

\maketitle
\begin{abstract}
This paper proposes and evaluates a Multi-Reference Relaxation Enforced Neighborhood Search (MRENS) heuristic within the SCIP solver.
This study marks the first integration and evaluation of MRENS in a full-fledged MILP solver, specifically coupled with the recently-introduced Lagromory separator for generating multiple reference solutions.
Computational experiments on the MIPLIB 2017 benchmark set show that MRENS, with multiple reference solutions, improves the solver’s ability to find higher-quality feasible solutions compared to single-reference approaches.
This study highlights the potential of multi-reference heuristics in enhancing primal heuristics in MILP solvers.

\keywords{Mixed-integer optimization \and Heuristics}

\end{abstract}

\section{Introduction}\label{sec:intro}
\renewcommand*{\thefootnote}{\fnsymbol{footnote}}
\footnotetext[1]{Corresponding author.}
\renewcommand*{\thefootnote}{\arabic{footnote}}

Primal heuristics aim to find feasible solutions to optimization problems at a lower computational cost than their exact counterparts but without any guarantees.
For mixed-integer linear programming (MILP), two main use cases spur the development of primal heuristics.
First, as standalone methods, heuristics allow practitioners to obtain feasible solutions in shorter amounts of time.
Second, as components of exact solution algorithms, heuristics provide feasible solutions for various purposes in these algorithms, e.g., for pruning nodes in a branch-and-cut algorithm.

In this work, we consider general MILP problems, which take the form:
\begin{align}
    \underset{{x}}{\mathrm{min}}\{{c}^{\intercal}{x} \;\; | \;\; {A}{x} \geq {b}, \;\; {l} \leq {x} \leq {u}, \;\; {x} \in \mathbb{Z}^{|\I|} \times \mathbb{R}^{n - |\I|} \}. \label{eq:milp}
\end{align}

Heuristics in MILP solvers can broadly be classified into four categories: rounding, diving, objective diving, and improvement heuristics.
Interested readers may refer to~\cite{berthold2006primal,achterberg2007constraint} for an overview.
The first three categories require a fractional solution as a reference solution, and the last category requires an incumbent feasible solution.
The fractional solution is often obtained from the linear programming (LP) relaxation of~\eqref{eq:milp}, which is defined as follows:
\begin{align}
    \underset{{x}}{\mathrm{min}}\{{c}^{\intercal}{x} \;\; | \;\; {A}{x} \geq {b}, \;\; {l} \leq {x} \leq {u}, \;\; {x} \in \mathbb{R}^{n} \}. \label{eq:lpr}
\end{align}
The usage of this fractional solution varies depending on the heuristic, e.g., \textit{relaxation enforced neighborhood search} (RENS)~\cite{berthold2006primal,achterberg2007constraint,berthold2014rens} and \textit{feasibility pump} (FP)~\cite{fischetti2005feasibility,achterberg2007improving,fischetti2009feasibility}.
Despite the existence of a wide variety of heuristics, there is one common aspect among them, i.e., all of them require and use a \textit{single} reference solution for their execution.

The literature for MILP heuristics exploiting \textit{multiple} reference solutions is sparse.
To the best of our knowledge, the only such work is~\cite{mexi2023using} in which the authors improved the FP heuristic by considering multiple reference solutions.
They proposed a new heuristic called mRENS and programmed it in a standalone implementation of the FP heuristic.
They observed performance improvements on various MILP test sets.
However, the idea was never tested and integrated into a full-fledged MILP solver.

The lack of literature using multiple reference solutions in a heuristic setting is also reflected in the optimization software landscape, where we are aware of no open-source implementations of such heuristics in an MILP solver.
Motivated by this, along with the recent integration of a new relax-and-cut framework-based separator in SCIP (more details in Section~\ref{sec:main}), we implement a new primal heuristic that considers multiple reference solutions within the state-of-the-art solver SCIP~\cite{bolusani2024scip}.
We call our new heuristic MRENS, similar to the heuristic in~\cite{mexi2023using}.

The differences between our proposed heuristic and the mRENS of~\cite{mexi2023using} are threefold.
First, we consider the multiple fractional solutions generated as a byproduct of the recently integrated \textit{Lagromory} separator~\cite{bolusani2024scip} in SCIP, whereas~\cite{mexi2023using} considers the multiple solutions generated in the pumping loop of the FP heuristic.
Second, the sub-MILPs of our heuristics differ in their feasible regions.
Finally, we use the working limits of the RENS heuristic in SCIP~\cite{achterberg2007constraint}, which are different from those of~\cite{mexi2023using}.

In the rest of this paper, we present the details of our heuristic in Section~\ref{sec:main}, our computational experiments and analysis in Section~\ref{sec:comp-results}, and our conclusions in Section~\ref{sec:conclusions}.

\section{Multi-Reference Relaxation Enforced Neighborhood Search (MRENS)}\label{sec:main}
RENS~\cite{berthold2006primal,achterberg2007constraint,berthold2014rens} is a rounding primal heuristic for MILPs that searches for an (integral) primal feasible solution of \eqref{eq:milp} in the neighborhood of a given single fractional solution, e.g., an optimal solution of \eqref{eq:lpr}.
MRENS, \textit{multi-reference relaxation enforced neighborhood search}, is a generalization of RENS in the sense that MRENS considers multiple fractional solutions to define the neighborhood.
We will now describe MRENS in Section~\ref{sec:mrens} and the generation procedure of multiple reference solutions in Section~\ref{sec:sol-gen}.

\subsection{MRENS}\label{sec:mrens}
Let $\fracx{0}$ be a fractional optimal solution of~\eqref{eq:lpr}.
A trivial rounding heuristic applied at $\fracx{0}$ may fail to find a primal feasible solution of~\eqref{eq:milp} most of the time because it ignores the linear constraints $A x \geq b$.
RENS instead solves the sub-MILP~\eqref{eq:milp-rens} and finds the best possible rounding of $\fracx{0}$ in the feasible region of~\eqref{eq:milp}.
\begin{align}
    \underset{{x}}{\mathrm{min}}\{{c}^{\intercal}{x} \; | \; {A}{x} \geq {b}, \; {l} \leq {x} \leq {u}, \; \lfloor \fracx{0}_j \rfloor    \leq x_j \leq \lceil\fracx{0}_j \rceil \; \forall j \in \I, \; {x} \in \mathbb{Z}^{|\I|} \times \mathbb{R}^{n - |\I|} \}. \label{eq:milp-rens}
\end{align}
While constructing~\eqref{eq:milp-rens}, the integer variables with integral values in $\fracx{0}$ are fixed to these values, and the domains of remaining integer variables are changed based on their fractional values in $\fracx{0}$.
Consequently, the sub-MILP~\eqref{eq:milp-rens} can be as computationally hard as the original MILP~\eqref{eq:milp}, contingent on the variable fixings and domain changes.
For example, if the original problem is a pure binary MILP, the sub-MILP constructed with a fully fractional reference solution results in the same original problem.
In practice, however, difficult sub-MILPs occur infrequently, and it has been observed empirically~\cite{berthold2014rens} that RENS typically produces over-restricted and thus infeasible sub-MILPs.

MRENS, by design, aims to overcome the observed issues of over-restricted sub-MILPs through the use of multiple reference solutions $\left\{\fracx{0}, \fracx{1}, \fracx{2}, \hdots, \fracx{k}\right\}$.
Let $\J = \{0, 1, 2, \hdots, k\}$, ${{x}_j}_{\textrm{min}} = \min\limits_{i \in \J} \fracx{i}_j$, and ${{x}_j}_{\textrm{max}} = \max\limits_{i \in \J} \fracx{i}_j$.
Then, in MRENS, we solve the sub-MILP~\eqref{eq:milp-mrens} instead of the sub-MILP~\eqref{eq:milp-rens}.
\begin{equation}
\begin{aligned}
    \underset{{x}}{\mathrm{min}}\{{c}^{\intercal}{x} \; | \; &{A}{x} \geq {b}, \; {l} \leq {x} \leq {u}, \; {x} \in \mathbb{Z}^{|\I|} \times \mathbb{R}^{n - |\I|},\\
										&\left\lceil {x_j}_{\textrm{min}} \right\rceil \leq x_j \leq \left\lfloor {x_j}_{\textrm{max}} \right\rfloor \; \forall j \in \I \;\; \textrm{if}\; {x_j}_{\textrm{max}} - {x_j}_{\textrm{min}} \geq 1.0,\\
										&\left\lfloor {x_j}_{\textrm{min}} \right\rfloor \leq x_j \leq \left\lceil {x_j}_{\textrm{max}} \right\rceil \; \forall j \in \I \;\; \textrm{if}\; {x_j}_{\textrm{max}} - {x_j}_{\textrm{min}} < 1.0 \}.
\end{aligned}
\label{eq:milp-mrens}
\end{equation}
The sub-MILPs~\eqref{eq:milp-mrens} and~\eqref{eq:milp-rens} are equivalent if we consider only the single fractional solution $\fracx{0}$, and otherwise, \eqref{eq:milp-mrens} is a relaxation of \eqref{eq:milp-rens}.
Accordingly,~\eqref{eq:milp-mrens} may be computationally more difficult to solve than~\eqref{eq:milp-rens} but also has a better chance of producing primal feasible solutions of~\eqref{eq:milp}.

\subsection{Generation of Multiple Reference Solutions}\label{sec:sol-gen}
Recently, the Lagromory separator (cutting plane generation routine) based on the relax-and-cut framework~\cite{fischetti2011relax} was implemented in SCIP 9.0~\cite{bolusani2024scip}.
The following steps describe one iteration of the execution process of this separator whenever it is called at a fractional optimal solution $\fracx{0}$ of~\eqref{eq:lpr}:
\begin{enumerate}[nosep]
\item Generate Gomory mixed-integer (GMI) cuts that separate $\fracx{0}$.
\item Add the cuts to the objective function of~\eqref{eq:lpr} in a Lagrangian fashion.
\item Update the Lagrangian multipliers and modify the objective of the LP.
\item Solve this new LP to obtain its optimal solution $\fracx{k}$ at iteration $k$.
\item If integral, save $\fracx{k}$ as a feasible solution to \eqref{eq:milp}. Otherwise, go to step 1.
\end{enumerate}
\noindent
Multiple such iterations are performed in a single call to the separator while keeping the feasible region of~\eqref{eq:lpr} intact.
This outputs a set of GMI cuts from multiple LP feasible bases of~\eqref{eq:lpr}.
Since it solves multiple LPs, this separator is computationally costly compared to other separators present in SCIP.
This computational cost motivated the search for ways to exploit different information generated throughout the separation algorithm, particularly on the primal side.
Thus, we propose to use these bases as the multiple reference solutions.

Using too many reference solutions can cause the sub-MILP to have a large feasible region inducing a prohibitive solving time.
Therefore, we only consider a maximum of three reference solutions, which follows the empirical results of \cite{mexi2023using}.
The reference solutions we select consist of the first solution $\fracx{0}$ and the last two solutions $\fracx{k-1}$ and $\fracx{k}$ from the sequence of solutions obtained by the Lagromory separator.

\section{Computational Results}\label{sec:comp-results}
We implemented and tested our heuristic in SCIP 9.0.0\footnote{SCIP commit hash: dadaf6a544b39ee20a64d1dde942b2b8b1164b7e}.
We did our experiments on a cluster equipped with Intel Xeon Gold 5122 CPUs and a limit of 96GB of RAM.
We used the MIPLIB 2017 benchmark library~\cite{MIPLIB2017} as our test set.
To mitigate the effects of performance variability~\cite{Lodi_2013}, we solved each problem with five different seeds, namely $\{0,1,2,3,4\}$, resulting in a total of 1200 problem-seed combinations, which we refer to as instances.

We solved each instance with a time limit of two hours and a memory limit of 50 GB.
We use the same working limits as RENS in SCIP, i.e., 5000 solving node limit and 500 stalling node limit (the maximum number of nodes an MILP solver can process without improving the incumbent solution of the sub-MILP).
We also require a minimum of 50\% integer variable fixings for heuristic execution while creating the sub-MILP and a minimum of 25\% total variable fixings after presolving the sub-MILP.
We modified SCIP to deactivate the built-in RENS heuristic and use the Lagromory separator as a routine for generating multiple reference solutions.

In the following, we compare two settings, ``MRENS'' and ``RENS''.
In the MRENS setting, we use three reference solutions as detailed in Section~\ref{sec:sol-gen}.
In the RENS setting, we use only the first solution $\fracx{0}$, which is the optimal solution of~\eqref{eq:lpr}.
Thus, the only difference between the RENS setting and the built-in RENS heuristic of SCIP is the frequency of execution, i.e., the RENS setting is only executed when the Lagromory separator generates a solution (refer to~\cite{bolusani2024scip} for additional details).

Table \ref{tab:1} provides a comparison of the two settings.
As expected, MRENS fixes fewer integer variables while constructing its sub-MILP~\eqref{eq:milp-mrens} compared to RENS.
Accordingly, MRENS is executed less often than RENS because of the requirement of minimum 50\% integer variable fixings.
MRENS is also more successful in finding feasible solutions than RENS.
More importantly, MRENS finds the best-known solution more frequently than RENS, with success rates of 17.0\% and 12.3\%, respectively.
\begin{table}[h]
\setlength{\tabcolsep}{3pt}
\centering
\begin{footnotesize}
\caption{Aggregated results comparing MRENS and RENS.
The columns in order of appearance: the heuristic setting used; the total number of heuristic calls over all instances; the percentages of heuristic calls where the heuristic was executed (i.e., the calls where at least 50\% of the integer variables were fixed in the sub-MILP), successfully found a solution, and found a new best solution; and the percentage of fixed integer variables in the sub-MILP averaged over all heuristic calls.}
\label{tab:1}
\begin{tabular}{lcccccc}
\toprule
\multicolumn{1}{c}{setting} & \# calls & \begin{tabular}[c]{@{}c@{}}\% of calls\\ heur. executed\end{tabular} &
\begin{tabular}[c]{@{}c@{}}\% of calls\\ solution found\end{tabular} & \begin{tabular}[c]{@{}c@{}}\% of calls\\ best found\end{tabular} 
& \begin{tabular}[c]{@{}c@{}}avg.~\% of fixed\\ int.~variables\end{tabular} \\
\midrule
MRENS	& 1883  & 78.5\%	& 30.5\%	& 17.0\%	& 73\%	\\
RENS	& 1848	& 83.0\%	& 24.1\%	& 12.3\%	& 77\%	\\
\bottomrule
\end{tabular}
\end{footnotesize}
\end{table}

In Table~\ref{tab:2}, we compare the overall performance of SCIP with MRENS and RENS settings.
Both settings solve almost the same number of instances within the time limit.
For the instances that were solved by both settings, MRENS is faster by 5\% and generates branch-and-bound trees with 3\% fewer nodes.
On the \textit{affected} instances where MRENS has impacted the solving process, MRENS outperforms RENS significantly both in terms of solution time and tree size.
Specifically, MRENS is 23\% faster and generates trees with 15\% fewer nodes.

\begin{table}[h]
\setlength{\tabcolsep}{4pt}
\centering
\begin{footnotesize}
\caption{Aggregated results comparing SCIP's performance with MRENS and RENS settings for the categories: all instances (``all''), instances solved by both settings (``both-solved''), instances where the solving path is different (``affected''), and affected instances solved by both settings (``affected-solved''). The columns refer to the number of instances in each category, the number of solved instances for each setting, the geometric mean of the runtime and nodes for each setting, and the relative quotients of those.}
\label{tab:2}
\begin{tabular}{lrrrrrrrrr}
\toprule
\multirow{2}{*}{subset} & \multirow{2}{*}{instances} & \multicolumn{3}{c}{MRENS} & \multicolumn{3}{c}{RENS} & \multicolumn{2}{c}{relative} \\
\cmidrule(lr){3-5} \cmidrule(lr){6-8} \cmidrule(lr){9-10}
			&		& \multicolumn{1}{c}{solved}	& \multicolumn{1}{c}{time}	& \multicolumn{1}{c}{nodes}	& \multicolumn{1}{c}{solved}	& \multicolumn{1}{c}{time}	& \multicolumn{1}{c}{nodes}	&	\multicolumn{1}{c}{time}	& \multicolumn{1}{c}{nodes} \\
\midrule
all			& 1192	& 632	& 1136	& -		& 633	& 1168	& -		& 0.97	& -		\\
both-solved	& 625	& 625	& 213	& 2608	& 625	& 224	& 2698	& 0.95	& 0.97	\\
\midrule
affected		& 142	& 134	& 322	& -		& 135	& 411	& -		& 0.78	& -		\\
affected-solved	& 127	& 127	& 228	& 5500	& 127	& 297	& 6475	& 0.77	& 0.85	\\
\bottomrule
\end{tabular}
\end{footnotesize}
\end{table}

Additionally, both settings are fast within the working limits; on average, they require less than one second of execution time per instance without accounting for the generation time of reference solutions.
For the instances that were solved to optimality, MRENS found an optimal solution for 38 out of 616, whereas RENS found an optimal solution for 13 out of 617.
For the instances that were not solved to optimality, MRENS found the best solution for 24 out of 562, whereas RENS found the best solution for 15 out of 559.

\section{Conclusion}\label{sec:conclusions}
Our results demonstrate that the MRENS setting finds both feasible and best-known solutions for MILPs more often than the RENS setting without much additional computational cost, indicating that using multiple reference solutions in the RENS framework is beneficial.
The concept may be extended to other heuristics, such as RINS~\cite{danna2005exploring}, or with other sources of reference solutions.
\section*{Acknowledgements}

The work for this article has been conducted in the Research Campus MODAL funded by the German Federal Ministry of Education and Research (BMBF) (grant number 05M14ZAM).

\bibliographystyle{splncs04}
\bibliography{refs}

\begin{thebibliography}{10}
\providecommand{\url}[1]{\texttt{#1}}
\providecommand{\urlprefix}{URL }
\providecommand{\doi}[1]{https://doi.org/#1}

\bibitem{achterberg2007constraint}
Achterberg, T.: Constraint integer programming. Ph.D. thesis, TU Berlin (2007)

\bibitem{achterberg2007improving}
Achterberg, T., Berthold, T.: Improving the feasibility pump. Discrete
  Optimization  \textbf{4}(1),  77--86 (2007)

\bibitem{berthold2006primal}
Berthold, T.: Primal heuristics for mixed integer programs. Master's thesis,
  Technische Universit{\"a}t Berlin (2006)

\bibitem{berthold2014rens}
Berthold, T.: {RENS}: the optimal rounding. Mathematical Programming
  Computation  \textbf{6},  33--54 (2014)

\bibitem{bolusani2024scip}
Bolusani, S., Besan{\c{c}}on, M., Bestuzheva, K., Chmiela, A., Dion{\'\i}sio,
  J., Donkiewicz, T., van Doornmalen, J., Eifler, L., Ghannam, M., Gleixner,
  A., et~al.: The {SCIP} optimization suite 9.0. arXiv preprint
  arXiv:2402.17702  (2024)

\bibitem{danna2005exploring}
Danna, E., Rothberg, E., Pape, C.L.: Exploring relaxation induced neighborhoods
  to improve {MIP} solutions. Mathematical Programming  \textbf{102},  71--90
  (2005)

\bibitem{fischetti2005feasibility}
Fischetti, M., Glover, F., Lodi, A.: The feasibility pump. Mathematical
  Programming  \textbf{104}(1),  91--104 (2005)

\bibitem{fischetti2009feasibility}
Fischetti, M., Salvagnin, D.: Feasibility pump 2.0. Mathematical Programming
  Computation  \textbf{1}(2),  201--222 (2009)

\bibitem{fischetti2011relax}
Fischetti, M., Salvagnin, D.: A relax-and-cut framework for {Gomory}
  mixed-integer cuts. Mathematical Programming Computation  \textbf{3},
  79--102 (2011)

\bibitem{MIPLIB2017}
Gleixner, A., Hendel, G., Gamrath, G., Achterberg, T., Bastubbe, M., Berthold,
  T., Christophel, P., Jarck, K., Koch, T., Linderoth, J., et~al.: {MIPLIB}
  2017: data-driven compilation of the 6th mixed-integer programming library.
  Mathematical Programming Computation pp. 1--48 (2021).
  \doi{10.1007/s12532-020-00194-3}

\bibitem{Lodi_2013}
Lodi, A., Tramontani, A.: Performance variability in mixed-integer programming.
  In: Theory Driven by Influential Applications, pp. 1--12. {INFORMS} (2013).
  \doi{10.1287/educ.2013.0112}

\bibitem{mexi2023using}
Mexi, G., Berthold, T., Salvagnin, D.: Using multiple reference vectors and
  objective scaling in the feasibility pump. EURO Journal on Computational
  Optimization  \textbf{11},  100066 (2023)

\end{thebibliography}

\end{document}